\documentclass[12pt]{amsproc}
\title{Divergent Ces\`aro and Riesz means of Jacobi and Laguerre expansions}
\author{Christopher Meaney}

\address{Department of Mathematics, Macquarie University, North Ryde
NSW 2109, Australia}

\email{chrism@maths.mq.edu.au}

\newtheorem{thm}{Theorem}[section]

\newtheorem{lem}[thm]{Lemma}

\theoremstyle{definition}

\theoremstyle{remark}


\begin{document}

 \begin{abstract} We show that for $\delta$ below certain critical indices there
are functions whose Jacobi or Laguerre expansions have  almost everywhere
divergent Ces\`aro and Riesz means of order $\delta$.
 \end{abstract} \maketitle

\section{Introduction}
\subsection{Orthogonal Expansions} Suppose that $(X,\mu)$ is a positive measure
space, $(\varphi_{n})_{n=0}^{\infty}$ is a orthogonal subset of 
$L^{2}(X,\mu)$, and  
$ {h_{n}=\|\varphi_{n}\|_{2}^{2}}$ for  all  $n\ge 0$.  
If $f$ is a function on $X$ for which all the products
$f.\overline{\varphi_{n} }$ are $\mu$-integrable then   $f$ has an 
orthogonal expansion
\begin{equation}\label{orthexp}
\sum_{n=0}^{\infty}
c_{n}(f) h_{n}^{-1}\varphi_{n}(x),
\end{equation}   with coefficients  
\begin{equation}\label{coeffs}
c_{n}(f)=\displaystyle\int_{X}f(x)\overline{\varphi_{n}(x)}\,d\mu(x),\quad\forall
n\ge 0.
\end{equation} 
\subsection{Ces\`aro  means}
As described in  Zygmund's book
\cite[pages 76--77]{MR38:4882},  the Ces\`aro  means of order $\delta$ of the
expansion (\ref{orthexp}) are defined by
\begin{equation}\label{cesexp}
\sigma_N^\delta f(x)=\sum_{n=0}^N
\dfrac{A^{\delta}_{N-n}}{A^{\delta}_N}c_{n}(f) h_{n}^{-1}\varphi_{n}(x),
\end{equation}
where $A^\delta_n=\displaystyle{{n+\delta}\choose n}$. Theorem 3.1.22 in
\cite{MR38:4882} says that if the Ces\`aro means converge, then the
terms of the series have controlled growth.

\begin{lem}\label{zygces} Suppose that   
$\displaystyle\lim_{N\to\infty}\sigma_N^\delta f(x)$ exists for some $x\in X$ and $\delta>
-1$. Then
$$
\left|  
c_{n}(f) h_{n}^{-1}\varphi_{n}(x)
\right| \le C_\delta\,  n^{\delta} \max_{0\le n\le N}\left|
\sigma_n^\delta f(x)
\right|,
\qquad \forall n\ge 0.
$$  
\end{lem} 
\subsection{Riesz means}
Hardy and Riesz \cite{HR} had proved a similar result for Riesz means. Recall
that the Riesz means of order $\delta\ge 0$ are defined for each  
$r>0$ by
\begin{equation}
\label{rieszexp}
S_{r}^{\delta}f(x)=\sum_{0\le n<r} 
\left( 1-\dfrac{k}{r}\right)^\delta 
c_{n}(f) h_{n}^{-1}\varphi_{n}(x).
\end{equation}
Theorem 21  of    \cite{HR}
tells us how the convergence of
$S_{r}^{\delta}f(x)$ controls the size of the partial sums 
$S_{r}^{0}f(x)$.
\begin{lem}\label{HardyRiesz}
Suppose that $f $ has an orthogonal expansion
 and for some  $\delta>0$ and $x\in X$ its Riesz means 
  $S_{r}^{\delta}f(x)$ converges to $c$
 as $r\to\infty$ then 
$$
\left| S_{r}^{0}f(x)-c\right|\le A_\delta\, r^{\delta} \sup_{0<t\le
r+1} \left| S_{t}^{\delta}f(x)\right|.
$$  
\end{lem}
Note that we can write
$$
c_{n}(f) h_{n}^{-1}\varphi_{n}(x)=\left(S_n^0f(x)-c\right)-\left(S_{n-1}^0
f(x)-c\right)={\bf O}(n^\delta)
$$
and obtain the same growth estimates as in Lemma \ref{zygces}.

Gergen\cite{0017.01103} wrote formulae relating the Riesz and Ces\`aro means of
order $\delta\ge 0$, from which the equivalence of the two methods of summation
follows.

\subsection{Uniform Boundedness}\label{ubp}
Suppose there is a $1<q\le \infty$ for which 
$\varphi_n\in L^q(X,\mu)$ for all
$n$. In addition, suppose that there is some
positive number $\lambda$ with 
\[
\| \varphi_n\|_q\ge c n^\lambda, \qquad\forall
n\ge 1.
\]
The formation of the coefficient  $f\mapsto c_n(f)$ is then a bounded linear
functional on the dual of $L^q(X,\mu)$ with norm bounded below by a constant
multiple of $n^\lambda$.
 The uniform boundedness principle implies that
for $p$ conjugate to $q$   and
each
$0\le
\varepsilon<\lambda$ there is an
$f\in L^p(X,\mu)$ so that 
 \begin{equation}\label{ubpinf}
c_n(f)/n^\varepsilon \to \infty \text{ as }
n\to\infty.
\end{equation}
 
\subsection{Cantor-Lebesgue Theorem}\label{CaLeTh}
The following argument is based on  
\cite[Section IX.1]{MR38:4882}.
Suppose we have a sequence of functions $F_n$ on an interval in the real line with
the asymptotic property  
$$
F_n(\theta)=c_n\left(\cos(M_n\theta +\gamma_n)+{\bf o}(1)\right),\qquad\forall n
\ge0
$$ 
uniformly on a set $E$ of finite positive measure, and 
with
$M_n\to  \infty$ as $n\to\infty$. Integrating $|F_n|^2$ over $E$ gives
\begin{eqnarray*}
 \int_E |F_n(\theta)|^2\,d\theta
 =|c_n|^2\left(\int_E \cos^2(M_n\theta+\gamma_n)\,d\theta
+{\bf o}(1)\right) \qquad\qquad \\
 =|c_n|^2\left(\dfrac{|E|}2 +\dfrac{e^{2i\gamma_n}}4
{\widehat{\chi_E}}(2M_n)+ \dfrac{e^{-2i\gamma_n}}4
{\widehat{\chi_E}}(-2M_n)+{\bf o}(1)\right).
\end{eqnarray*}
The Riemann-Lebesgue Theorem  \cite[Thm. II.4.4]{MR38:4882}   says
that the Fourier transforms
 ${\widehat{\chi_E}}(\pm 2M_n)\to 0$ as 
$M_n\to\infty$.
If we know that there is some function $G$ for which  $| F_n(\theta)| \le G(n)$
uniformly on
$E$ for all
$n$ then there is an $n_0>0$ for which 
$$
\dfrac{ |E|}4 |c_n|^2 \le  \int_E |F_n(\theta)|^2\,d\theta \le
G(n)^2|E|,\quad\forall n\ge n_0.
$$
This shows that $|c_n|\le 2G(n)$ for all $n\ge n_0$.

\section{Jacobi Polynomials}
\subsection{Notation}
 Fix real numbers $\alpha\ge\beta\ge-1/2$, with $\alpha>-1/2$,
and let
$\mu$ denote the measure on
$[-1,1]$ defined by 
$$ d\mu(x)=(1-x)^\alpha (1+x)^\beta
\,dx.$$ Let $P^{(\alpha,\beta)}_n(x)$ be the Jacobi polynomial of
degree $n$ associated to the pair $(\alpha,\beta)$ as in  Szeg{\H{o}}'s book
\cite{MR46:9631}. 
Then $\left( P^{(\alpha,\beta)}_n\right)_{n=0}^\infty$ is an orthogonal subset
of $L^2([-1,1],\mu)$.
Equation $(4.3.3)$ in \cite{MR46:9631} shows
that the normalization terms
${h_{n}^{(\alpha,\beta)}=  \left\|
P^{(\alpha,\beta)}_n
\right\|^2_2 }$ satisfy
\begin{equation}\label{jach}
h^{(\alpha,\beta)}_{n}\sim 
c_{\alpha,\beta}n^{-1} \text{ as }n\to \infty.
\end{equation}
The Jacobi polynomial expansion of $f\in
L^1(\mu)$ is
\[
\sum_{n=0}^\infty  c_n(f) \left(
h_{n}^{(\alpha,\beta)}\right)^{-1} P^{(\alpha,\beta)}_n(x),
\]
 with
  coefficients $c_n(f) =\int_{-1}^{1}f(x) P^{(\alpha,\beta)}_n(x)\,d\mu(x)$. We
take $\alpha$ and $\beta$ as fixed and use $\sigma^\delta_N f(x)$ and
$S_r^\delta f(x)$ to denote the Ces\`aro and Riesz means of this expansion,
respectively.
 
\subsection{Asymptotics} Theorem 8.21.8 in Szeg{\H{o}}'s
book\cite{MR46:9631} gives the following asymptotic behaviour for the
Jacobi polynomials.
\begin{lem}\label{darboux} For $\alpha\ge \beta\ge -1/2$ and
$\varepsilon>0$ the following estimate holds uniformly for all 
${\varepsilon\le \theta\le \pi-\varepsilon}$ and $n\ge 1$.
\begin{equation}\label{jacas}
P^{(\alpha,\beta)}_n(\cos\theta)=n^{-1/2}k(\theta) \cos\left(M_n
\theta +\gamma\right)+{\bf O}\left( n^{-3/2}
\right).
\end{equation}
Here $k(\theta)=\pi^{-1/2}\left(\sin(\theta/2)
\right)^{-\alpha-1/2} \left(\cos(\theta/2)
\right)^{-\beta-1/2}$,\hfill\break  ${M_n=n+(\alpha+\beta+1)/2}$, and
${\gamma=-(\alpha+1/2)\pi/2}$
\end{lem}  From
Egoroff's theorem and Lemma \ref{zygces} we can say that if 
$ \sigma_{N}^\delta f(x)$  converges on a set of positive measure in
$[-1,1]$ then there is a set of positive measure  $E$ on which 
\begin{equation}\label{main1}
\left| c_{n} n^{(1/2)-\delta}\left( \cos\left( M_n\theta
+\gamma\right) +{\bf O}(n^{-1 })\right)
\right|\le A
\end{equation} uniformly for $\cos\theta\in E$. The argument of subsection
\ref{CaLeTh}  shows that
\begin{equation}\label{main2}
\left|  c_n n^{(1/2)-\delta}\right|\le A, \qquad\forall  n\ge 1.
\end{equation} 
\subsection{Norm estimates}
Next we recall the calculation of Lebesgue norms of
Jacobi polynomials, according to  Markett \cite{MR85d:42025} and  Dreseler and
Soardi
\cite{MR84d:42029}. Equation (2.2) in \cite{MR85d:42025} gives the following
lower bounds on these norms.
\begin{lem}\label{Markett} For real numbers $\alpha\ge\beta\ge-1/2$,
with $\alpha>-1/2$, ${1\le q<\infty}$, and ${r>-1/q}$,
$$
\left(
\int_{0}^1
\left| P_{n}^{(\alpha,\beta)}(x) \, (1-x)^r
\right|^q \, dx
\right)^{1/q}\sim 
\begin{cases} 
\displaystyle n^{-1/2} & \text{ if }\quad  r>\alpha/2+1/4-1/q,\\ 
\displaystyle n^{-1/2} \left(\log n\right)^{1/q} &  
\text{ if }\quad  r=\alpha/2+1/4-1/q,\\ 
\displaystyle n^{\alpha-2r-2/q} & \text{ if }\quad 
r<\alpha/2+1/4-1/q.\\
\end{cases}
$$
\end{lem}

 \subsection{Main Result} 
There are critical indices, as used in
\cite{MR84c:42040}, 
$$ p_c=\dfrac{4(\alpha+1)}{(2\alpha+3)}\quad
\text{and its conjugate  }p_c'=\dfrac{4(\alpha+1)}{(2\alpha+1)}.
$$
Taking $r=\alpha/q$ in Lemma \ref{Markett} we have that
\begin{equation}\label{jaclpnorm}
\left\| P_{n}^{(\alpha,\beta)} \right\|_{L^q(\mu)}>
\left(
\int_{0}^1
\left| P_{n}^{(\alpha,\beta)}(x)\right|^q \, (1-x)^\alpha
 \, dx
\right)^{1/q} \sim n^{\alpha-2\alpha/q-2/q}
\end{equation} for $\alpha/q<\alpha/2+1/4-1/q$. This last inequality
can be rewritten as
\begin{equation} q>\dfrac{4(\alpha+1)}{2\alpha+1}=p_c'.
\end{equation}

We can now prove that below the critical index there are functions
with almost everywhere divergent Ces\`aro and  Riesz means.
\begin{thm} For real numbers $\alpha\ge\beta\ge-1/2$, with
$\alpha>-1/2$, 
$$ 1\le p <p_c=\dfrac{4(\alpha+1)}{(2\alpha+3)},
\text{ and  }
  0\le \delta <{\dfrac{ \left(  2\alpha+2 \right) }
    { p}} - {\dfrac{\left(2\alpha+3
         \right)   }{2}},
$$ there is an $f\in L^p(\mu)$, supported in $[0,1]$, whose Ces\`aro
means $\sigma_N^\delta f(x)$ and Riesz means $S_r^\delta f(x)$ are divergent
almost everywhere on
$[-1,1]$.
\end{thm}

\proof  Let $ q $ be conjugate to $p$, so that $1/p=(q-1)/q$.
Suppose that 
\[
\delta <{\dfrac{ \left(  2\alpha+2 \right) }
    { p}} - {\dfrac{\left(2\alpha+3
         \right)   }{2}} = {\dfrac{ \left(  2\alpha+2 \right)\left( q-1\right) }
    { q}} - {\dfrac{\left(2\alpha+3
         \right)   }{2}}=\alpha+\dfrac12-\dfrac{(2\alpha+2)}q,
\] then
$$
\delta
-\dfrac12< \alpha -\dfrac{(2\alpha+2)}q=\alpha-\dfrac{2\alpha}q-\dfrac2q,
$$
which is the exponent of $n$ in the inequality
(\ref{jaclpnorm}). Now apply the argument given in
subsection~\ref{ubp}. The norms of the Jacobi
polynomials in Lemma \ref{Markett} are calculated over
$[0,1]$ and so we can find $f$ in
$L^p([-1,1],\mu)$, supported on $[0,1]$, for
which the coefficients satisfy
$$
c_n(f)/n^{\delta-1/2}\to\infty,\quad\text{ as }
n\to \infty.
$$
Combine this with Lemmas \ref{zygces} and \ref{HardyRiesz} and the
argument around inequality (\ref{main2}) to see that for this $f$ both 
$ 
\sigma_N^\delta f(x)  $ and $  S_r^\delta f(x)$ 
are divergent almost everywhere. This argument follows  the
methods used in \cite{MR84c:42040,MR89i:42030,MR90b:42034,MR99j:43003}. \qed

\subsection{Remarks}
Convergence results above the critical index are contained in the work of
Bonami and Clerc \cite{MR49:3461}, Colzani, Taibleson and Weiss
\cite{MR86g:43012}, and Chanillo and Muckenhoupt\cite{MR93g:42018}.  In
particular, in \cite[Thm. 1.4]{MR93g:42018} it is shown that for 
$$ 1\le p <p_c=\dfrac{4(\alpha+1)}{(2\alpha+3)},
\text{ and  }
   \delta ={\dfrac{ \left(  2\alpha+2 \right) }
    { p}} - {\dfrac{\left(2\alpha+3
         \right)   }{2}},
$$
the maximal operator $f\mapsto \sup_{N\ge 0}\left| \sigma_{N}^\delta
f(x)\right|$ is of weak type $(p,p)$.

\section{Laguerre Functions}
\subsection{Notation}
For each $\alpha>-1$ let $\mu_\alpha$ be the
measure on
$[0,\infty)$ defined by 
$$
d\mu_\alpha(x)=e^{-x}x^\alpha \, dx.
$$
We denote by $L^{(\alpha)}_n(x)$ the Laguerre
polynomial of degree $n$, as in  
\cite[Chpt. 5]{MR46:9631}. The $L^2(\mu_\alpha)$-norms of these satisfy the
identity
$$
h_n^{(\alpha)}=\|
L_n^{(\alpha)}\|_{L^2(\mu_\alpha)}^2=\Gamma(\alpha+1){n+\alpha
\choose n} \sim n^{\alpha}.
$$
 Fej\'er's formula  
\cite[Thm. 8.22.1]{MR46:9631} gives the asymptotic properties of
these polynomials. For each $\alpha>-1$ and
$0<\varepsilon<\omega$,
\begin{equation}\label{lagas}
L^{(\alpha)}_n(x)=\dfrac{e^{x/2}}{\pi^{ 1/2} x^{
\alpha/2}}
n^{\alpha/2-1/4}\cos\left(2(nx)^{1/2}-\alpha\pi/2-\pi/4
\right)
+{\bf O}\left(
n^{\alpha/2-3/4}
\right),
\end{equation}
uniformly in $x\in[\varepsilon,\omega]$.
 The corresponding normalized functions
are
 \begin{equation}\label{LaFun}
{\mathcal L}_n^\alpha (x)
=\sqrt{\dfrac{\Gamma(n+1)}{\Gamma(n+\alpha+1)}}
\,e^{-x/2}
x^{a/2}L_n^{(\alpha)}(x),\qquad\forall x\ge 0,
n\ge 0.
\end{equation}
These provide an orthonormal subset of
$L^2([0,\infty))$, where the half line carries
  Lebesgue measure.

\subsection{Norm Estimates}Markett\cite[Lemma 1]{MR83j:40004} has calculated the
Lebesgue norms of the Laguerre functions, for $\alpha>-1/2$
\begin{equation}\label{LagFNorm}
\left\| 
 {\mathcal L}_n^\alpha 
\right\|_q\sim 
\begin{cases}
n^{ 1/q-1/2},& \forall 1\le q <2,\\
n^{-1/4}(\log n)^{1/4},& \text{if } q=4,\\
n^{-1/q} ,& \forall 4<q\le \infty.
\end{cases}
\end{equation}

\subsection{Divergence Result}

\begin{thm}
If $\alpha>-1/2$, $p>4$
and $0<\delta<1/4-1/p$ then there is a function
$f\in L^p(0,\infty)$
whose Laguerre expansion $$\sum_{n=0}^\infty c_n(f){\mathcal L}_n^\alpha(x)$$ has
Ces\`aro and Riesz means of order $\delta$ which diverge almost
everywhere. 

\end{thm}
\proof
Suppose that the expansion
$ 
\sum_{n=0}^\infty c_n(f) {\mathcal L}_n^\alpha(x)
$ 
is either Ces\`aro or Riesz summable of order
$\delta$ on a set of positive measure in
$[0,\infty)$. Then Lemma  \ref{zygces} or Lemma \ref{HardyRiesz} implies that
\begin{equation}\label{couldbe}
c_n(f) {\mathcal L}_n^\alpha(x) ={\bf O}(n^\delta) 
\end{equation}
on a set of positive measure.
When equations (\ref{lagas}) and ({\ref{LaFun}}) are combined with the argument
of subsection
\ref{CaLeTh} we find that
\begin{equation}\label{delbnd}
c_n(f)  ={\bf O}(n^{\delta+1/4}).
\end{equation}
The case when $\delta=0$  is Lemma 2.3 in Stempak's paper \cite{MR2001g:42057}.
Suppose that 
$$
\dfrac1q-\dfrac12>\delta+\dfrac14,
$$
so that
$\delta<\dfrac1q-\dfrac34=\dfrac{4-3q}q$. If
$\dfrac1q=1-\dfrac1p$ then this inequality is 
$ 
\delta<\dfrac14-\dfrac1p$. 
The argument  of subsection \ref{ubp} shows that if $p>4$
and $\delta<1/4-1/p$ then there is a function
$f\in L^p(0,\infty)$ for which the inequality (\ref{delbnd}) fails, 
$$
c_n(f)/n^{\delta+1/4}\to\infty \text{ as } n\to \infty.
$$
The Laguerre expansion  of this function has Ces\`aro and Riesz
means of order $\delta$ which diverge almost
everywhere. \qed

\subsection{Remarks}
There is an extensive treatment of almost everywhere convergence results for
Laguerre expansions in \cite{MR2001m:42051}. In particular, \cite[Thm.
1.20]{MR2001m:42051} implies that if $p>4$ and $\delta\ge 1/4-1/p$ then all
$f\in L^p( 0,\infty )$ have almost everywhere convergent Ces\`aro means  of order
$\delta$.


\bibliographystyle{amsplain}
\providecommand{\bysame}{\leavevmode\hbox to3em{\hrulefill}\thinspace}

\end{document}